\documentclass[12pt]{amsart}

\usepackage{amsmath,xspace,amssymb,mathrsfs}
\usepackage{color}

\input xy
\xyoption{all}
\xyoption{2cell}
\UseAllTwocells
\CompileMatrices

\newcommand{\Spec}{\operatorname{Spec}}
\renewcommand{\phi}{\varphi}

\newcommand{\Ker}{\operatorname{Ker}}

\newcommand{\Max}{\operatorname{Max}}

\newcommand{\Ann}{\operatorname{Ann}}
\newcommand{\Hom}{\operatorname{Hom}}

\newcommand{\Supp}{\operatorname{Supp}}

\newcommand{\tr}{\operatorname{tr}}

\newtheorem{proposition}{Proposition}[section]
\newtheorem{lemma}[proposition]{Lemma}

\newtheorem{corollary}[proposition]{Corollary}
\newtheorem{theorem}[proposition]{Theorem}

\theoremstyle{definition}

\newtheorem{remark}[proposition]{Remark}

\begin{document}

\title[Pure ideals and trace ideals]{Some results on pure ideals and trace ideals of projective modules}

\author[A. Tarizadeh]{Abolfazl Tarizadeh}
\address{Department of Mathematics, Faculty of Basic Sciences, University of Maragheh \\
P. O. Box 55136-553, Maragheh, Iran.
 }
\email{ebulfez1978@gmail.com}

\date{}

\begin{abstract} Let $R$ be a commutative ring with the unit element. It is shown that an ideal $I$ in $R$ is pure if and only if $\Ann(f)+I=R$ for all $f\in I$. If $J$ is the trace of a projective $R$-module $M$, we prove that $J$ is generated by the ``coordinates" of $M$ and $JM = M$. These lead to a few new results and alternative proofs for some known results.
\end{abstract}

\maketitle

\section{Introduction and Preliminaries}

The concept of the trace ideals of modules has been the subject of research by some mathematicians around late 50's until late 70's and has again been active in recent years (see, e.g. \cite{Beckwith}, \cite{Dao et al.}, \cite{Herbera}, \cite{Herzog}, \cite{Jondrup}, \cite{lindo}, \cite{Vasconcelos 2} and \cite{Whitehead}). This paper
deals with some results on the trace ideals of projective modules. We begin with a few results on pure ideals which are used in their comparison with trace ideals in the sequel. After a few preliminaries in the present section, in section 2 a new characterization of pure ideals is given (Theorem \ref{Remark 030}) which
is followed by some corollaries. Section 3 is devoted to the trace ideal of
projective modules. Theorem \ref{Remark 0501} gives a characterization of the trace ideal of a projective module in terms of the ideal generated by the ``coordinates" of the elements of the module. This characterization enables us to deduce some new results on the trace ideal of projective modules like the statement on the trace ideal of the tensor product of two modules for which one of them is projective (Corollary \ref{Corollary VI89}), and some alternative proofs for a few known results such as Corollary \ref{Corollary IV} which shows that the trace ideal of a projective module is a pure ideal. \\
In this paper all rings are assumed to be commutative with the unit element. \\

The following theorem will be used in a number of our results.

\begin{theorem}\label{lemma 8} Every finitely generated flat module over a local ring is free.
\end{theorem}

{\bf Proof.} See \cite[Tag 00NZ]{Johan} or \cite[Theorem 7.10]{Matsumura}. $\Box$ \\

\begin{remark}\label{Remark I} Let $F$ be a free $R-$module with a basis $(e_{k})_{k\in K}$ and $M$ an $R-$module. If $x\in F\otimes_{R}M$ then there exists a unique sequence $(m_{k})\in\bigoplus\limits_{k\in K}M$ such that $x=\sum\limits_{k}e_{k}\otimes m_{k}$. In fact, the map $\bigoplus\limits_{k\in K}M\rightarrow F\otimes_{R}M$ given by $(m_{k})\rightsquigarrow\sum\limits_{k}e_{k}\otimes m_{k}$ is an isomorphism of $R-$modules. This, in particular, implies that if $\phi:R\rightarrow S$ is a morphism of rings then $(e_{k}\otimes1)_{k\in K}$ is a basis for the free $S-$module $F\otimes_{R}S$.
\end{remark}

A projective $R-$module is also called $R-$projective. The same terminology is used for free and flat modules as well. \\

\section{Pure ideals}

An ideal $I$ of a ring $R$ is called a pure ideal if the canonical ring map
$R\rightarrow R/I$ is a flat ring map. Pure ideals are studied in commutative and non-commutative algebra. In this section we give some results on pure ideals. \\

Theorem \ref{Remark 030} translates the $R$-flatness of $R/I$ in terms of coprime-ness of the ideal $I$ with the annihilator of any of its elements.
The main motivations comes from \cite[Theorem 7.13]{Matsumura} and first properties of an absolutely flat ring, a ring $R$ where every $R-$module is $R-$flat (see Corollary \ref{Corollary II}(iii)).

\begin{theorem}\label{Remark 030} An ideal $I$ of a ring $R$ is a pure ideal if and only if $\Ann(f)+I=R$ for all $f\in I$.
\end{theorem}

{\bf Proof.} First assume $R/I$ is $R-$flat. Suppose there is some $f\in I$ such that $\Ann(f)+I\neq R$. Thus there exists a prime ideal $\mathfrak{p}$ of $R$ such that $\Ann(f)+I\subseteq\mathfrak{p}$. Therefore, by Theorem \ref{lemma 8}, $I_{\mathfrak{p}}=0$. So there exists some $s\in R\setminus\mathfrak{p}$ such that $sf=0$. But this is a contradiction. \\
Conversely, let $\phi:M\rightarrow N$ be an injective morphism of $R-$modules. To prove the assertion it suffices to show that the induced map $M/IM\rightarrow N/IN$ given by $m+IM\rightsquigarrow \phi(m)+IN$ is injective. If $\phi(m)\in IN$ then we may write $\phi(m)=\sum\limits_{i=1}^{n}f_{i}x_{i}$ where $f_{i}\in I$ and $x_{i}\in N$ for all $i$. By the hypothesis, there are elements $b_{i}\in\Ann(f_{i})$ and $c_{i}\in I$ such that $1=b_{i}+c_{i}$. It follows that $1=(b_{1}+c_{1})(b_{2}+c_{2})\cdots(b_{n}+c_{n})=b+c$ where $b=b_{1}b_{2}\cdots b_{n}$ and $c\in I$. Thus $\phi(m)=b\phi(m)+c\phi(m)=\phi(cm)$. Therefore $m=cm\in IM$.  $\Box$

\begin{remark} Here we give an elementary proof (without using of  Theorem \ref{lemma 8}) for the implication ``$\Rightarrow$" of Theorem \ref{Remark 030}. If $R/I$ is $R-$flat, then the exact sequence of $R-$modules $\xymatrix{0\ar[r]&I\ar[r]&R}$ gives rise to the exact sequence: $$\xymatrix{0\ar[r]&I/I^{2}\simeq I\otimes_{R}R/I\ar[r]&R\otimes_{R}R/I\simeq R/I.}$$
Note that the above morphism sends each pure tensor $a\otimes(r+I)\in I\otimes_{R}R/I$ into $ra+I=0$. Hence, its image is zero and so $I=I^{2}$. Now if $f\in I$ then the exact sequence $\xymatrix{0\ar[r]&Rf\ar[r]&I}$ gives us the exact sequence: $$\xymatrix{0\ar[r]&Rf\otimes_{R}R/I\ar[r]&I\otimes_{R}R/I\simeq I/I^{2}=0.}$$
Thus $R/(\Ann(f)+I)\simeq R/\Ann(f)\otimes_{R}R/I\simeq Rf\otimes_{R}R/I=0$ and so $\Ann(f)+I=R$. \\
\end{remark}

\begin{corollary}\label{Corollary I} The annihilator of a finitely generated flat module is a pure ideal.
\end{corollary}

{\bf Proof.} Let $M$ be a finitely generated flat module over a ring $R$ with the annihilator $I=\Ann_{R}(M)$. Suppose there is some $f\in I$ such that $\Ann(f)+I\neq R$. Thus there exists a prime ideal $\mathfrak{p}$ of $R$ such that $\Ann(f)+I\subseteq\mathfrak{p}$. Clearly $M_{\mathfrak{p}}$ is a nonzero finitely generated flat module over the local ring $R_{\mathfrak{p}}$. Then, by Theorem \ref{lemma 8}, $M_{\mathfrak{p}}$ is a free $R_{\mathfrak{p}}-$module. So $I_{\mathfrak{p}}=\Ann_{R_{\mathfrak{p}}}(M_{\mathfrak{p}})=0$. Hence there is some $s\in R\setminus\mathfrak{p}$ such that $sf=0$. Thus $s\in\Ann(f)\subseteq\mathfrak{p}$ which is a contradiction. Therefore by Theorem \ref{Remark 030}, $I$ is a pure ideal. $\Box$ \\

In view of Theorem \ref{Remark 030} and Corollary \ref{Corollary I}, the well-known characterization of the pure ideals could be further extended (see \cite[Tag 04PS]{Johan}).

\begin{corollary}\label{Corollary II} Let $I$ be an ideal of a ring $R$. Then the following statements are equivalent. \\
$\mathbf{(i)}$ $I$ is a pure ideal. \\
$\mathbf{(ii)}$ $I=\{f\in R: \Ann(f)+I=R\}$. \\
$\mathbf{(iii)}$ If $f\in I$ then there exists some $g\in I$ such that $f(1-g)=0$. \\
$\mathbf{(iv)}$  $\Supp(I)=\Spec(R)\setminus V(I)$. \\
$\mathbf{(v)}$ If $\mathfrak{p}$ is a prime ideal of $R$, then either $I_{\mathfrak{p}}=0$ or $I_{\mathfrak{p}}=R_{\mathfrak{p}}$.
\end{corollary}

{\bf Proof.} $\mathbf{(i)}\Rightarrow\mathbf{(iv)}:$ If $I_{\mathfrak{p}}\neq0$ then there exists some $f\in I$ such that $f/1\neq0$. This yields that $\Ann(f)\cap(R\setminus\mathfrak{p})=\emptyset$. It follows that $\Ann(f)\subseteq\mathfrak{p}$. By Theorem \ref{Remark 030}, $\Ann(f)+I=R$. Therefore $\mathfrak{p}\in\Spec(R)\setminus V(I)$. \\
$\mathbf{(iv)}\Rightarrow\mathbf{(v)}:$ Easy. \\
$\mathbf{(v)}\Rightarrow\mathbf{(i)}:$ If $\Ann(f)+I\neq R$ for some $f\in I$ then there exists a prime ideal $\mathfrak{p}$ of $R$ such that $\Ann(f)+I\subseteq\mathfrak{p}$. This yields that $f/1\neq0$ and so $I_{\mathfrak{p}}=R_{\mathfrak{p}}$. This means that $I\nsubseteq\mathfrak{p}$, a contradiction. \\
The remaining implications follow easily by applying Theorem \ref{Remark 030}. $\Box$ \\

Motivated by Theorem \ref{Remark 030}, then we call an ideal $I$ of a ring $R$ \emph{strongly pure} if $\Ann(f)+Rf=R$ for all $f\in I$. \\

Every regular ideal (i.e., generated by a set of idempotents) is a pure ideal, see \cite[Lemma 8.4]{A. Tarizadeh}. But the converse does not hold, see \cite[Example 4.7]{A. Tarizadeh 2}. However, regarding strongly pure ideals we have the following result.

\begin{proposition} Every strongly pure ideal is a regular ideal.
\end{proposition}

{\bf Proof.} Let $I$ be a strongly pure ideal of a ring $R$. If $f\in I$ then by hypothesis, there exist $g\in Rf$ and $h\in\Ann(f)$ such that $g+h=1$. Then $g(1-g)=gh=0$ and so $g$ is an idempotent. We also have $f=fg$. Thus $Rf=Rg$ is a regular ideal. Therefore $I$ is a regular ideal, because
$I=\sum\limits_{f\in I}Rf$. $\Box$ \\

The converse of above result does not hold. In fact, if $e$ is an idempotent of a ring $R$, then $Re$ is not necessarily a strongly pure ideal. As a specific example, the regular ideal $\mathbb{Z}$ is not strongly pure, since $\Ann(2)+2\mathbb{Z}=2\mathbb{Z}\neq\mathbb{Z}$. \\

In summary, in a given ring we have the following inclusions of sets: \\
Strongly pure ideals $\subseteq$ Regular ideals $\subseteq$ Pure ideals. \\

By \cite[Theorem 2.2(v)]{A. Tarizadeh 2}, a ring
$R$ is zero-dimensional if and only if for each $f\in R$ there exists a natural number $n\geqslant1$ such that $\Ann(f^{n})+Rf=R$.
Also note that if $f$ is a member of a ring $R$, then $\Ann(f)\subseteq\Ann(f^{2})\subseteq\cdots$.
Motivated by these observations, then we may define an ideal $I$ of a ring $R$ a \emph{quasi-pure} ideal if for each $f\in I$ there exists a natural number $n\geqslant1$ such that $\Ann(f^{n})+I=R$. Clearly an ideal $I$ is quasi-pure if and only if for each $f\in I$ there exists some $g\in I$ such that $f(1-g)$ is nilpotent. Similarly above, we call an ideal $I$ of a ring $R$  \emph{strongly quasi-pure} if for each $f\in I$ there exists a natural number $n\geqslant1$ such that $\Ann(f^{n})+Rf=R$. In a reduced ring $R$, these new notions are respectively coincide with the ``pure ideal'' and ``strongly pure ideal'' notions, because $\Ann(f)=\Ann(f^{n})$ for all $f\in R$ and $n\geqslant2$. If $R$ is a zero-dimensional ring, then each ideal of $R$ is strongly quasi-pure.

\begin{lemma}\label{Lemma iv quasi-pure} If each maximal ideal of a ring $R$ is quasi-pure, then $R$ is zero-dimensional.
\end{lemma}

{\bf Proof.} If $\mathfrak{p}$ is a prime ideal of $R$, then $\mathfrak{p}\subseteq\mathfrak{m}$ for some
maximal ideal $\mathfrak{m}$ of $R$. If the inclusion is strict, then we may choose some $f\in\mathfrak{m}$ such that $f\notin\mathfrak{p}$. By hypothesis, there exists some $g\in\mathfrak{m}$ such that $f(1-g)$ is nilpotent. If follows that $1-g\in\mathfrak{p}$ which is a contradiction. Hence, every prime ideal of $R$ is maximal. $\Box$ \\

If $R$ is an absolutely flat (i.e., reduced zero-dimensional) ring, then each ideal of $R$ is strongly pure, because it is well known that $R$ is absolutely flat if and only if it is von-Neumann regular ring (i.e., each $f\in R$ can be written as $f=f^{2}g$ for some $g\in R$). Here ``reducedness'' is crucial. For example, in the zero-dimensional ring $R=\mathbb{Z}/8\mathbb{Z}$, the ideal $(2)=\{0,2,4,6\}$ is not pure and so it is not strongly pure. The following result shows that the reverse also holds, even under the much weaker condition.

\begin{proposition} If each maximal ideal of a ring $R$ is a pure ideal, then $R$ is absolutely flat.
\end{proposition}

{\bf Proof.} By Lemma \ref{Lemma iv quasi-pure}, $R$ is zero-dimensional. It remains to show that $R$ is reduced. If $f\in R$ is nilpotent, then $f\in\mathfrak{m}$ for all $\mathfrak{m}\in\Max(R)$. By hypothesis, $\mathfrak{m}=\Ker\pi_{\mathfrak{m}}$ where $\pi_{\mathfrak{m}}:R\rightarrow R_{\mathfrak{m}}$ is the canonical ring map. So there exists some $g_{\mathfrak{m}}\in R\setminus\mathfrak{m}$ such that $fg_{\mathfrak{m}}=0$. Clearly the ideal $\big(g_{\mathfrak{m}}: \mathfrak{m}\in\Max(R)\big)$ is the whole ring $R$. Thus we may write $1=\sum\limits_{i=1}^{n}r_{i}g_{i}$ where $g_{i}:=g_{\mathfrak{m}_{i}}$ for all $i$. It follows that $f=\sum\limits_{i=1}^{n}r_{i}fg_{i}=0$. $\Box$ \\

For further information and applications of the pure ideals we refer the interested reader to the literature, especially to the recent works \cite{A. Tarizadeh} and \cite{A. Tarizadeh 2}-\cite{A. Tarizadeh 4}. \\

\section{The trace ideal of a projective module}

In this section we provide some results present our contributions on the
trace ideal of a projective module. \\

Recall that the \emph{trace ideal} of an $R-$module $M$ is the ideal $$\tr_{R}(M):=\sum\limits_{f\in M^{\ast}}f(M)$$
where $M^{\ast}=\Hom_{R}(M,R)$. \\

It is easy to see that the trace ideal of a free module $F$ is either the zero ideal or the whole ring, according to whether $F=0$ or $F\neq0$. \\

Let $M$ be a projective $R-$module and let $F$ be a free $R-$module which admits $M$ as a direct summand. Then there exists an $R-$submodule $N$ of $F$ such that $F=M+N$ and $M\cap N=0$. Let $(e_{k})_{k\in K}$ be a basis of $F$. Then each $m\in M$ can be written uniquely as $m=\sum\limits_{k}r_{m,k}e_{k}$ where $r_{m,k}=0$ for all but a finite number of indices $k$. The scalars $r_{m,k}$ are called the coordinates of $m$ with respect to the above basis. Note that the coordinates of each element of $M$ may be changed by passing from the basis $(e_{k})$ to another basis for $F$. But the following result  (motivated by \cite[p. 132, Proposition 3.1]{Cartan-Eilenberg}) shows that the ideal generated by the coordinates of the members of $M$ is independent of the choice of a free module $F$ and a basis for it. In fact, it is the trace ideal of $M$.

\begin{theorem}\label{Remark 0501} Let $M$ be a projective $R-$module. Then $J=\tr_{R}(M)$ is generated by the coordinates of the elements of $M$. In particular, $JM=M$.
\end{theorem}

{\bf Proof.} Let $J'$ be the ideal of $R$ generated by the coordinates of the elements of $M$. Using the above setup, then for each $i\in K$, by the universal property of free modules, there exists a (unique) morphism of $R-$modules $g_{i}:F\rightarrow R$ such that $g_{i}(e_{k})=\delta_{i,k}$ where $\delta_{i,k}$ is the Kronecker delta. Setting $f_{i}:=g_{i}\circ j$ where $j:M\rightarrow F$ is the canonical injection. Then $r_{m,i}=f_{i}(m)\in J$ for all $m\in M$ and all $i\in K$. So $J'\subseteq J$. Also each $m=\sum\limits_{k}r_{m,k}x_{k}$ where each $e_{k}=x_{k}+y_{k}$ with $x_{k}\in M$ and $y_{k}\in N$. In particular, we have $JM=M$. If $\phi:M\rightarrow R$ is a morphism of $R-$modules then $\phi(x_{k})\in R$ for all $k\in K$, and so $\phi(m)\in J'$.
Therefore $J=J'$. $\Box$ \\

In the sequel we will observe that the above theorem has some interesting consequences. \\

The following result is well known, see \cite[Propositions A.1. and A. 3. and Theorem A.2. ]{Auslander-Goldman} and \cite[after Proposition 1.3]{Vasconcelos}. We give an alternative proof which seems shorter.

\begin{corollary}\label{lemma 1} The annihilator of a finitely generated projective module is generated by an idempotent element.
\end{corollary}

{\bf Proof.} If $M$ is a projective $R-$module then by Theorem \ref{Remark 0501}, $JM=M$ where $J=\tr_{R}(M)$. If moreover, $M$ is a finitely generated $R-$module then by the so called ``the determinant trick", there exists some $f\in J$ such that $1-f\in I:=\Ann(M)$. Using the expression of the generators of $M$ as a linear combinations of
the basis elements $e_{k}$ considered in the proof of Theorem \ref{Remark 0501}, it follows that $IJ=0$. Setting $g:=1-f$. Then we have $g^{2}=g(1-f)=g$. Thus $g$ is an idempotent. If $h\in I$ then $hg=h(1-f)=h$ and so $h\in Rg$. Therefore $I=Rg$. $\Box$

\begin{corollary}\label{Corollary VII980} Let $M$ be an $R-$module and $x\in M$. Then $Rx$ is $R-$projective if and only if $\Ann(x)$ is generated by an idempotent of $R$.
\end{corollary}

{\bf Proof.} The implication ``$\Rightarrow$'' is deduced from Corollary \ref{lemma 1}. Conversely, by hypothesis there exists an idempotent $e\in R$ such that $\Ann(x)=Re$. We have the canonical isomorphisms of $R-$modules $Rx\simeq R/Re\simeq R(1-e)$ and $R(1-e)$ is a direct summand of $R$. Hence, $Rx$ is $R-$projective. $\Box$

\begin{theorem}\label{Corollary III} If $\phi:R\rightarrow S$ is morphism of rings and $M$ is a projective $R-$module, then $\tr_{S}(M\otimes_{R}S)=\tr_{R}(M)S$.
\end{theorem}

{\bf Proof.} Clearly $\tr_{R}(M)S\subseteq\tr_{S}(M\otimes_{R}S)$. To see the reverse inclusion, using the notations of Theorem \ref{Remark 0501}, then $M\otimes_{R}S$ is a direct summand of the free $S-$module $F\otimes_{R}S$. But it is well known that the family $e_{k}\otimes1$ is a basis for this free module, see Remark \ref{Remark I} or \cite[Corollary in page 26]{Northcott}. By Theorem \ref{Remark 0501}, $\tr_{S}(M\otimes_{R}S)$ is generated by the coordinates of the pure tensors $m\otimes s$ of $M\otimes_{R}S$. We have $m\otimes s=\sum\limits_{k}\phi(r_{m,k})s(e_{k}\otimes 1)$ and each $\phi(r_{m,k})s\in\tr_{R}(M)S$. Therefore $\tr_{S}(M\otimes_{R}S)=\tr_{R}(M)S$. $\Box$ \\

The following result is well known, a sketch of its proof can be found in \cite[p. 16]{Jondrup} also see \cite[p. 269]{Vasconcelos 2}. We provide an alternative proof.

\begin{corollary}\label{Corollary IV} The trace ideal of a projective module over a commutative ring is a pure ideal.
\end{corollary}

{\bf Proof.} If $M$ is a projective $R-$module then for each prime ideal $\mathfrak{p}$ of $R$, by Theorem \ref{Corollary III}, $\tr_{R_{\mathfrak{p}}}(M_{\mathfrak{p}})=J_{\mathfrak{p}}$ where $J=\tr_{R}(M)$. But it is well known that every projective module over a local ring is free, see \cite[Theorem 2]{Kaplansky}. This yields that either $J_{\mathfrak{p}}=0$ or $J_{\mathfrak{p}}=R_{\mathfrak{p}}$. Hence, by Corollary \ref{Corollary II}, $J$ is a pure ideal. $\Box$ \\

It is important to notice that the ``commutativity'' assumption of Corollary \ref{Corollary IV} is crucial. If we drop this assumption then the assertion does not hold anymore, see \cite[Example 1.2]{Jondrup}. \\

In \cite[Proposition 1.1]{Jondrup} the converse of Corollary \ref{Corollary IV} is also proved which states that if $I$ is a pure ideal of a ring $R$ then there exists a projective $R-$module whose trace ideal is equal to $I$. \\

\begin{corollary}\label{Corollary VI89} If $M$ is a projective $R-$module then for each $R-$module $N$, $\tr_{R}(M\otimes_{R}N)=JJ'=J\cap J'$ where $J=\tr_{R}(M)$ and $J'=\tr_{R}(N)$.
\end{corollary}

{\bf Proof.} It is easy to see that $JJ'\subseteq\tr_{R}(M\otimes_{R}N)\subseteq J\cap J'$ for every two $R-$modules $M$ and $N$. If $M$ is a projective $R-$module then by Corollary \ref{Corollary IV},
$J$ is a pure ideal and so $JJ'=J\cap J'$. $\Box$

\begin{lemma}\label{Lemma I0100} Let $I$ be an ideal of a ring $R$ such that $\tr_{R}(I)=R$. Then $I$ is a finitely generated ideal.
\end{lemma}

{\bf Proof.} We may write $1=\sum\limits_{i=1}^{n}f_{i}(x_{i})$ where $f_{i}\in\Hom_{R}(I,R)$ and $x_{i}\in I$ for all $i$. If $y\in I$ then $y=\sum\limits_{i=1}^{n}f_{i}(y)x_{i}$. Hence, $I=(x_{1},\ldots,x_{n})$ is a finitely generated ideal. $\Box$ \\

Finally, we provide an alternative proof to the following result.

\begin{corollary}\cite[Proposition 1.1]{Vasconcelos 2}\label{Corollary V} Let $I$ be an ideal of a ring $R$ which is not contained in any minimal prime ideal of $R$. If $I$ is a projective $R-$module, then it is a finitely generated ideal.
\end{corollary}

{\bf Proof.} By Lemma \ref{Lemma I0100}, it suffices to show that $\tr_{R}(I)=R$. If not, then there exists a maximal ideal $\mathfrak{m}$ of $R$ such that $J:=\tr_{R}(I)\subseteq\mathfrak{m}$. There exists a minimal prime ideal $\mathfrak{p}$ of $R$ such that $\mathfrak{p}\subseteq\mathfrak{m}$. By the hypotheses, there exists some $x\in I$ such that $x\notin\mathfrak{p}$. We have $I=IJ\subseteq J$, see Theorem \ref{Remark 0501}. By Corollary \ref{Corollary IV}, $J$ is a pure ideal and so by Theorem \ref{Remark 030}, $\Ann(x)+J=R$. But $\Ann(x)\subseteq\mathfrak{p}$ and so $\Ann(x)+J\subseteq\mathfrak{m}$.
This is a contradiction and we win. $\Box$ \\

\textbf{Acknowledgements.} The author would like to give heartfelt thanks to the referee for very careful reading of the paper and for his/her excellent comments and modifications which significantly improved the paper.

\end{document}